\documentclass[a4paper, reqno]{amsart}
\pdfoutput=1
\title[{}]{Atiyah--Bott formula and homological block of Seifert fibered homology $3$-spheres}

\author{Shoma Sugimoto}
\address{Faculty of Mathematics, Kyushu University, Fukuoka 819-0386, Japan}
\email{shomasugimoto361@gmail.com}

\usepackage{latexsym}
\usepackage{amsmath}
\usepackage{amsfonts}
\usepackage{amssymb}
\usepackage{amsxtra}
\usepackage{mathrsfs}
\usepackage{eucal}
\usepackage{tikz} \usetikzlibrary{calc} \tikzset{>=latex} \usetikzlibrary{backgrounds}
\usepackage[all]{xy}
\usepackage{color}
\usepackage{braket}
\usepackage{graphics, setspace}
\usepackage{etoolbox, textcomp}
\usepackage{comment}
\usepackage{changepage}
\usepackage{xcolor}
\definecolor{rouge}{rgb}{0.85,0.1,.4}
\definecolor{bleu}{rgb}{0.1,0.2,0.9}
\definecolor{violet}{rgb}{0.7,0,0.8}
\usepackage[colorlinks=true,linkcolor=bleu,urlcolor=violet,citecolor=rouge]{hyperref}
\usepackage{cleveref}
\usepackage{lscape}

\usepackage{tikz}			
\usetikzlibrary{matrix}
\usetikzlibrary{automata, arrows.meta, positioning,cd}
\usepackage{dynkin-diagrams}
\usetikzlibrary{cd}			

\newtheorem{definition}{Definition}[section]

\newtheorem{remark}[definition]{Remark}

\numberwithin{equation}{section}

\newcommand{\Z}{\mathbb{Z}}
\newcommand{\R}{\mathbb{R}}
\newcommand{\C}{\mathbb{C}}

\newcommand{\im}{\operatorname{Im}}

\newcommand{\ch}{\operatorname{ch}}
\newcommand{\tr}{\operatorname{tr}}

\newcommand{\g}{\mathfrak{g}}

\newcommand{\sll}{\mathfrak{sl}}

\newcommand\doi[2]{\href{http://dx.doi.org/#1}{#2}}

\newcommand{\rk}{\mathrm{rank}\g}

\allowdisplaybreaks
\begin{document}
\maketitle
\begin{abstract}
We show that the homological block of Seifert fibered homology $3$-spheres are computed by the iterated use of Atiyah--Bott fixed point formula.
\end{abstract}

\section{introduction}\label{section introduction}
In \cite{GPPV}, Gukov--Pei--Putrov--Vafa introduced $q$-series called homological block for any plumbed $3$-manifold.
Homological blocks have a remarkable property that their radial limits coincides with the Witten--Reshetikhin--Turaev invariant, as conjectured in \cite{GPPV}.
It had been proved in various cases \cite{AnM,FIMT,MMur,Mur1}, and a general proof was given very recently in \cite{Mur2}.
It is also known that ``spoiled'' modular forms such as false/mock theta functions appear as homological blocks, and their modular transformation has been studied \cite{MTer,BMM}. 

On the other hand, homological blocks are expected to appear as characters of modules of logarithmic vertex operator algebras.
In fact, it is known that \cite{CFGH,CFGHP} the homological blocks of $3$ or $4$-fibered Seifert $3$-manifolds are given by certain linear combinations of the irreducible modules of $(1,p)$-log VOA \cite{AdM,CRW,FT} and $(p,p')$-log VOA \cite{FGST,TW}, respectively.
The irreducible modules (in general, baby Verma modules) of $(1,p)$-log VOA have a geometric realization by Feigin--Tipunin \cite{FT,Sug,Sug2}.
They are given by the space of global sections $H^0(G\times_BV_{\lambda})$ of certain equivariant bundle $G\times_BV_{\lambda}$ over the flag variety $G/B$ whose fiber space is the corresponding irreducible module $V_{\lambda}$ of the lattice VOA $V_{\sqrt{p}Q}$.
The key point is that by combining the cohomology vanishing $H^{n>0}(G\times_BV_{\lambda})=0$, which is a natural consequence of the theory of Felder exactness \cite{Sug,CNS2}, and Atiyah--Bott fixed point formula \cite{AB}, the character of $H^0(G\times_BV_{\lambda})$ \cite{FT,Sug} is computed from that of $V_{\lambda}$.
Recently, this technique was also used to compute the characters of the irreducible modules of $(p,p')$-log VOA \cite{KSug} and of the inverse Hamiltonian reduction of $(1,p)$-log VOA \cite{CNS,CNS2}.

In this paper, we show that the homological blocks of $(N+2)$-fibered Seifert homology $3$-spheres are computed by repeating the technique above ($N\geq 1$).
In other words, iterated use of Atiyah--Bott formula to the corresponding lattice theta function gives the desired homological block.
This result suggests the existence of ``nested Feigin--Tipunin construction" of log VOAs and their modules of the form
\footnote{If taken literally, this is nonsense, but the point is that nested character computation can be performed using Atiyah--Bott formula.
See \S\ref{B-submodules} and \ref{character formula} below.}
\begin{align}
``H^0(G\times_BH^0(G\times_B\cdots H^0(G\times_BV_{\vec{\lambda}})))"
\end{align}
whose construction and the study of representation theory are for future works.

\subsection*{Acknowledgements}
The author wishes to express his gratitude to Yuji Terashima, Kazuhiro Hikami, Toshiki Matsusaka, Yuya Murakami, Tomoyuki Arakawa, Shigenori Nakatsuka for useful comments and discussions. 
The author is also deeply grateful to Sergei Gukov for his interest in this idea.
This work was supported by JSPS KAKENHI Grant number 22J00951.

\section{Nested Feigin--Tipunin construction}
\subsection{Notation}
Unless otherwise noted, the ground field is an complex field $\C$.
Let $\mathfrak{g}$ be a simply-laced simple Lie algebra, and $\mathfrak{g}=\mathfrak{n}_-\oplus\mathfrak{h}\oplus\mathfrak{n}_+$ the triangular decomposition, $\mathfrak{h}$ the Cartan subalgebra, $\mathfrak{b}=\mathfrak{n}_-\oplus\mathfrak{h}$ the Borel subalgebra, $G$ and $B$ the semisimple, simply-connected algebraic groups corresponding to $\mathfrak{g}$ and $\mathfrak{b}$, respectively. 
Set $I:=\{1,\dots,\mathrm{rank}\g\}$.
We adopt the standard numbering for the simple roots $\{\alpha_i\}_{i\in I}$ of $\mathfrak{g}$ and $\{\varpi_i\}_{i\in I}$ denotes the corresponding fundamental weights.
For $j\in I$, denote $\sll_2^j\simeq\sll_2$ (resp. $\mathrm{SL}_2^j$) the Lie subalgebra of $\g$ generated by $\{e_j,h_j,f_j\}$ (resp. corresponding subgroup). 
For $J\subseteq I$, $P_J$ means the corresponding parabolic subgroup.
Denote by $\Delta^{\pm}$ the sets of positive roots and negative roots of $\mathfrak{g}$, respectively.
Let $Q$ be the root lattice of $\mathfrak{g}$, $P$ the weight lattice of $\mathfrak{g}$, $P_+$ the set of dominant integral weights of $\mathfrak{g}$. 
Denote by $(\cdot , \cdot)$ the invariant form of $\mathfrak{g}$ normalized as $(\alpha_i,\alpha_i)=2$ for any $1\leq i\leq \rk$, $W$ the Weyl group of $\mathfrak{g}$ generated by the simple reflections $\{\sigma_i\}_{i\in I}$, $(c^{ij})$ the inverse matrix to the Cartan matrix of $\mathfrak{g}$, $\rho$ the half sum of positive roots, $\theta$ the highest root of $\mathfrak{g}$, $h$ the Coxeter number of $\mathfrak{g}$.
For $\sigma\in W$, $l(\sigma)$ denotes the length of $\sigma$.
For $\beta\in P_+$, $L(\beta)$ denote the irreducible $\mathfrak{g}$-module with highest weight $\beta$.
Let $\{h_i\}_{i\in I}$ be the basis of $\mathfrak{h}$ corresponding to the simple roots by $(\cdot , \cdot)$.
For $\mu\in P$, $\C(\mu)$ denotes the one-dimensional $\mathfrak{b}$-module such that the action of $\mathfrak{n}_-$ is trivial and $h_i$ acts as $(\alpha_i,\mu)\operatorname{id}$ for $i\in I$.
For a $B$-module $V$ and $\mu\in P$, we use the letter $V(\mu)$ for the $B$-module $V\otimes \C(\mu)$.
For $\beta\in\mathfrak{h}^\ast$, $V^{h=\beta}$ donotes the weight space of the $B$-module $V$.

\subsection{Lattice VOA}
Let us consider the case of simply-laced $\mathfrak{g}$.
Let $p_1,\dots,p_N$ be a family of coprime integers and $p:=p_1\cdots p_N$.
We consider the lattice VOA $V_{\sqrt{p}Q}$ associated with the rescaled root lattice $\sqrt{p}Q$ with the conformal vector
\begin{align*}
\omega:=\tfrac{1}{2}\sum_{1\leq i,j\leq \mathrm{rank}\g}c^{ij}{\alpha_i}_{(-1)}\alpha_j+\sqrt{p}Q_0\rho_{(-2)}\mathbf{1},\quad Q_0:=\tfrac{1}{p_1}-\tfrac{1}{p_2}-\cdots-\tfrac{1}{p_N},
\end{align*}
where $(c^{ij})_{1\leq i,j\leq \mathrm{rank}\g}$ is the inverse of Cartan matrix and $\rho$ is the Weyl vector, namely, the half sum of positive roots.
For $1\leq m\leq N$, set
\begin{align*}
&\Lambda_{p_m}=\{\lambda_m=\sum_{i=1}^{\mathrm{rank}\g}\tfrac{1-r_{m,i}}{p_m}\varpi_i~|~1\leq r_{m,i}\leq p_m\},\\
&\vec{\lambda}=(\lambda_1,\cdots,\lambda_N)
:=\sqrt{p}(\lambda_1-\lambda_{2}\cdots-\lambda_N)
\end{align*}
Then we have the complete set of irreducible $V_{\sqrt{p}Q}$-modules 
\begin{align*}
V_{\hat{\lambda};\vec{\lambda}}=V_{\hat{\lambda};\lambda_1,\dots,\lambda_N}:=V_{\sqrt{p}(Q-\hat\lambda)+\vec{\lambda}},
\quad \text{$\hat{\lambda}$ is a minuscle weight,}
\end{align*}
and the conformal weight $\Delta_{\sqrt{p}\beta+\vec{\lambda}}$ of $e^{\sqrt{p}\beta+\vec{\lambda}}\in V_{\hat{\lambda};\vec{\lambda}}$ is given by
\begin{align*}
\Delta_{\sqrt{p}\beta+\vec{\lambda}}=
\tfrac{p}{2}|\beta+\vec{\lambda}-Q_0\rho|^2.
\end{align*}
When $\hat\lambda$ is unimportant, $V_{\hat\lambda,\vec{\lambda}}$ is simply written as $V_{\vec{\lambda}}$.
For a vector space $M$ with the conformal grading, the character of $M$ is defined by $\ch_qM=\tr_{M}q^{L_0-\tfrac{c}{24}}$.

\subsection{Felder complex}
\begin{figure}[htbp]
\centering
\begin{tikzpicture}[scale=.4]

\node[draw] (-65r) at (10,0) {$5$};

\node[draw]
(-43r) at (8,4) {$3$};
\node[draw]
(-44r) at (8,2) {$4$};
\node[draw]
(-45r) at (8,0) {$5$};

\node[draw]
(-21r) at (6,8) {$1$};
\node[draw]
(-22r) at (6,6) {$2$};
\node[draw]
(-23r) at (6,4) {$3$};
\node[draw]
(-24r) at (6,2) {$4$};
\node[draw]
(-25r) at (6,0) {$5$};

\node[draw]
(00r) at (4,10) {$0$};
\node[draw]
(01r) at (4,8) {$1$};
\node[draw]
(02r) at (4,6) {$2$};
\node[draw]
(03r) at (4,4) {$3$};
\node[draw]
(04r) at (4,2) {$4$};
\node[draw]
(05r) at (4,0) {$5$};

\node[draw]
(22r) at (2,6) {$2$};
\node[draw]
(23r) at (2,4) {$3$};
\node[draw]
(24r) at (2,2) {$4$};
\node[draw]
(25r) at (2,0) {$5$};

\node[draw]
(44r) at (0,2) {$4$};
\node[draw]
(45r) at (0,0) {$5$};

\draw[arrows=->] (01r)--(-21r);

\draw[arrows=->] (22r)--(02r);
\draw[arrows=->] (02r)--(-22r);

\draw[arrows=->] (23r)--(03r);
\draw[arrows=->] (03r)--(-23r);
\draw[arrows=->] (-23r)--(-43r);

\draw[arrows=->] (44r)--(24r);
\draw[arrows=->] (24r)--(04r);
\draw[arrows=->] (04r)--(-24r);
\draw[arrows=->] (-24r)--(-44r);

\draw[arrows=->] (45r)--(25r);
\draw[arrows=->] (25r)--(05r);
\draw[arrows=->] (05r)--(-25r);
\draw[arrows=->] (-25r)--(-45r);
\draw[arrows=->] (-45r)--(-65r);

\node[draw]
(-54r) at (24,2) {$4$};
\node[draw]
(-55r) at (24,0) {$5$};

\node[draw]
(-32r) at (22,6) {$2$};
\node[draw]
(-33r) at (22,4) {$3$};
\node[draw]
(-34r) at (22,2) {$4$};
\node[draw]
(-35r) at (22,0) {$5$};

\node[draw]
(-10r) at (20,10) {$0$};
\node[draw]
(-11r) at (20,8) {$1$};
\node[draw]
(-12r) at (20,6) {$2$};
\node[draw]
(-13r) at (20,4) {$3$};
\node[draw]
(-14r) at (20,2) {$4$};
\node[draw]
(-15r) at (20,0) {$5$};

\node[draw]
(11r) at (18,8) {$1$};
\node[draw]
(12r) at (18,6) {$2$};
\node[draw]
(13r) at (18,4) {$3$};
\node[draw]
(14r) at (18,2) {$4$};
\node[draw]
(15r) at (18,0) {$5$};

\node[draw]
(33r) at (16,4) {$3$};
\node[draw]
(34r) at (16,2) {$4$};
\node[draw]
(35r) at (16,0) {$5$};

\node[draw]
(55r) at (14,0) {$5$};

\draw[arrows=->] (11r)--(-11r);

\draw[arrows=->] (12r)--(-12r);
\draw[arrows=->] (-12r)--(-32r);

\draw[arrows=->] (33r)--(13r);
\draw[arrows=->] (13r)--(-13r);
\draw[arrows=->] (-13r)--(-33r);

\draw[arrows=->] (34r)--(14r);
\draw[arrows=->] (14r)--(-14r);
\draw[arrows=->] (-14r)--(-34r);
\draw[arrows=->] (-34r)--(-54r);

\draw[arrows=->] (55r)--(35r);
\draw[arrows=->] (35r)--(15r);
\draw[arrows=->] (15r)--(-15r);
\draw[arrows=->] (-15r)--(-35r);
\draw[arrows=->] (-35r)--(-55r);

\node (a) at (10,5) {$\ $}; 
\node (x) at (14,5) {$\ $};
\node (b) at (10,4) {$\ $}; 
\node (y) at (14,4) {$\ $};
\node (c) at (12,6) {$Q^{[\lambda_m]}_j$}; 
\node (z) at (12,3) {$Q^{[\sigma_j\ast\lambda_m]}_j$};
\draw[arrows=->] (a)--(x);
\draw[arrows=->] (y)--(b);

\node (h) at (-2,-2) {$h_j=$};
\node (h) at (0,-2) {$4$};
\node (h) at (2,-2) {$2$};
\node (h) at (4,-2) {$0$};
\node (h) at (6,-2) {$-2$};
\node (h) at (8,-2) {$-4$};
\node (h) at (10,-2) {$-6$};
\node (h) at (14,-2) {$5$};
\node (h) at (16,-2) {$3$};
\node (h) at (18,-2) {$1$};
\node (h) at (20,-2) {$-1$};
\node (h) at (22,-2) {$-3$};
\node (h) at (24,-2) {$-5$};
\end{tikzpicture}
\caption{Felder exact sequence.}
\label{fig: felder complex}
\end{figure}
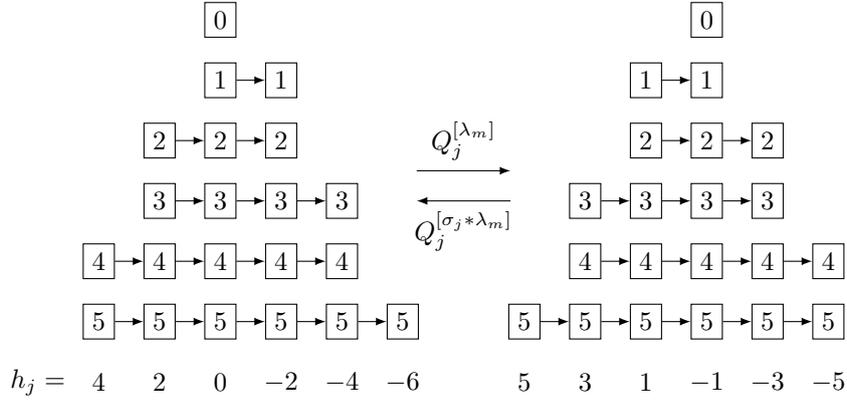

For the Weyl group $W$, $\sigma\in W$ and $\lambda_m\in\Lambda_m$, set
\begin{align*}
\sigma\ast\lambda_m=\sigma(\lambda_m+\tfrac{1}{p_m}\rho)-\tfrac{1}{p_m}\rho.
\end{align*}
For $\gamma\in\tfrac{1}{p_m}P$, denote $[\gamma]$ by the representative of $\gamma$ in $\Lambda_{p_m}\simeq \tfrac{1}{p_m}P/P$.
Then $\lambda_m\rightarrow[\sigma\ast\lambda_m]$ above defines a $W$-action on $\Lambda_m$.
We introduce a useful notation
\begin{align*}
\epsilon_{\lambda_m}(\sigma):=\tfrac{1}{p_m}(\sigma\ast\lambda_m-[\sigma\ast\lambda_m])\in P.
\end{align*}
For the properties of the $W$-action and $\epsilon_{\lambda_m}(\sigma)$, see \cite{Sug,CNS2}.
\begin{definition}
Let $\{V_{\lambda_m}\}_{\lambda_m\in\Lambda_m}$  
be a family of weight $B$-modules equipped with linear operators 
$Q^{[\lambda_m]}_j\colon V_{\lambda_m}\rightarrow V_{\sigma_j\ast\lambda_m}$ $(j\in I)$.
We call 
$\{V_{\lambda_m},Q^{[\lambda_m]}_j\}_{\lambda_m\in\Lambda_m,j\in I}$
Felder exact if
\begin{enumerate}
\item 
$Q^{[\lambda_m]}_j\colon V_{\lambda_m}\rightarrow V_{\sigma_j\ast\lambda_m}(\epsilon_{\lambda}(\sigma_j))$ is a $B$-module homomorphism.
\item 
$W_{j,\lambda_m}:=\ker_{V_{\lambda_m}}Q^{[\lambda_m]}_j$ is an $P_j$-submodule of $V_{\lambda_m}$.
\item 
For $j\in I$ and $\lambda_m\not\in\Lambda_m^{\sigma_j}$, we have the short exact sequence of $B$-modules
\begin{align*}
0\rightarrow
W_{j,\lambda_m}
\rightarrow
V_{\lambda_m}
\rightarrow
W_{j,\sigma_j\ast\lambda_m}(\epsilon_\lambda(\sigma_j))
\rightarrow 0,
\end{align*}
which is called Felder exact sequence.
\end{enumerate}
\end{definition}
A Felder exact sequence is illustrated as Figure \ref{fig: felder complex}.
Here, each
\begin{tikzpicture}
[baseline=(01l.base)]
\node[draw] (01l) at (0,0) {$k$};
\end{tikzpicture}
represents a family of weight vectors of simple $\mathrm{SL}_2^j$-module $L(k\varpi)$ (or $B$-module $L(k\varpi)(-\varpi)$) with highest weight $k\varpi$.

In the same manner as \cite{Sug}, we have
\begin{itemize}
\item 
The evaluation map $H^0(G\times_BV_{\lambda_m})\rightarrow V_{\lambda_m}$, $s\mapsto s(\mathrm{id}_{G/B})$ is an injective $B$-module homomorphism \cite[Theorem 4.14]{Sug}.
In particular, we can identify $H^0(G\times_BV_{\lambda_m})$ with the maximal $G$-submodule of $V_{\lambda_m}$.
\item 
The image of the evaluation map above is contained in $\bigcap_{j\in I}\ker Q^{[\lambda_m]}_j$.
We have $H^0(G\times_BV_{\lambda_m})\simeq \bigcap_{j\in I}\ker Q^{[\lambda_m]}_j$ iff $\lambda_m$ satisfies \cite[(100)]{Sug}.
\item 
(Borel--Weil--Bott type duality, \cite[Theorem 4.8]{Sug})
If $(-p_m\lambda_m+\rho,\theta)\leq p_m$, then we have a $G$-module isomorphism 
\begin{align*}
H^n(G\times_BV_{\lambda_m})
\simeq 
H^{n+l(w_0)}(G\times_BV_{w_0\ast\lambda_m}(-\rho)).
\end{align*}
In particular, $H^{n>0}(G\times_BV_{\lambda_m})=0$.
\end{itemize}
\begin{remark}
The Borel--Weil--Bott type duality above is an analog of \cite[Corollary 5.5(c)]{J}.
From the results in \cite{Sug,Sug2}, a strong parallel relationship between the theory of Felder complex and that of algebraic groups in finite characteristic (or that of quantum group at root of unity) is expected.
The author expects that the Borel--Weil--Bott type duality only holds for the case $(-p_m\lambda_m+\rho,\theta)\leq p_m$, but the cohomology vanishing $H^{n>0}(G\times_BV_{\lambda_m})=0$ holds for any $\lambda_m$, because it would be an analog of Kempf's vanishing theorem \cite[Proposition 4.5]{J}.
\end{remark}

\subsection{$B$-submodules}\label{B-submodules}
\begin{figure}[htbp]
\centering
\begin{tikzpicture}[scale=.4]

\node[fill=red!50] (-64r) at (10,2) {$4$};
\node[draw] (-65r) at (10,0) {$5$};

\node[fill=red!50] (-42r) at (8,6) {$2$};
\node[draw]
(-43r) at (8,4) {$3$};
\node[draw]
(-44r) at (8,2) {$4$};
\node[draw]
(-45r) at (8,0) {$5$};

\node[fill=red!50] (-20r) at (6,10) {$0$};
\node[draw]
(-21r) at (6,8) {$1$};
\node[draw]
(-22r) at (6,6) {$2$};
\node[draw]
(-23r) at (6,4) {$3$};
\node[draw]
(-24r) at (6,2) {$4$};
\node[draw]
(-25r) at (6,0) {$5$};

\node[draw]
(00r) at (4,10) {$0$};
\node[draw]
(01r) at (4,8) {$1$};
\node[draw]
(02r) at (4,6) {$2$};
\node[draw]
(03r) at (4,4) {$3$};
\node[draw]
(04r) at (4,2) {$4$};
\node[draw]
(05r) at (4,0) {$5$};

\node[draw]
(22r) at (2,6) {$2$};
\node[draw]
(23r) at (2,4) {$3$};
\node[draw]
(24r) at (2,2) {$4$};
\node[draw]
(25r) at (2,0) {$5$};

\node[draw]
(44r) at (0,2) {$4$};
\node[draw]
(45r) at (0,0) {$5$};

\draw[arrows=->] (01r)--(-21r);

\draw[arrows=->] (22r)--(02r);
\draw[arrows=->] (02r)--(-22r);

\draw[arrows=->] (23r)--(03r);
\draw[arrows=->] (03r)--(-23r);
\draw[arrows=->] (-23r)--(-43r);

\draw[arrows=->] (44r)--(24r);
\draw[arrows=->] (24r)--(04r);
\draw[arrows=->] (04r)--(-24r);
\draw[arrows=->] (-24r)--(-44r);

\draw[arrows=->] (45r)--(25r);
\draw[arrows=->] (25r)--(05r);
\draw[arrows=->] (05r)--(-25r);
\draw[arrows=->] (-25r)--(-45r);
\draw[arrows=->] (-45r)--(-65r);

\node[fill=red!50] (-70r) at (26,0) {$5$};

\node[fill=red!50] (-53r) at (24,4) {$3$};
\node[draw]
(-54r) at (24,2) {$4$};
\node[draw]
(-55r) at (24,0) {$5$};

\node[fill=red!50] (-31r) at (22,8) {$1$};
\node[draw]
(-32r) at (22,6) {$2$};
\node[draw]
(-33r) at (22,4) {$3$};
\node[draw]
(-34r) at (22,2) {$4$};
\node[draw]
(-35r) at (22,0) {$5$};

\node[draw]
(-10r) at (20,10) {$0$};
\node[draw]
(-11r) at (20,8) {$1$};
\node[draw]
(-12r) at (20,6) {$2$};
\node[draw]
(-13r) at (20,4) {$3$};
\node[draw]
(-14r) at (20,2) {$4$};
\node[draw]
(-15r) at (20,0) {$5$};

\node[draw]
(11r) at (18,8) {$1$};
\node[draw]
(12r) at (18,6) {$2$};
\node[draw]
(13r) at (18,4) {$3$};
\node[draw]
(14r) at (18,2) {$4$};
\node[draw]
(15r) at (18,0) {$5$};

\node[draw]
(33r) at (16,4) {$3$};
\node[draw]
(34r) at (16,2) {$4$};
\node[draw]
(35r) at (16,0) {$5$};

\node[draw]
(55r) at (14,0) {$5$};

\draw[arrows=->] (11r)--(-11r);

\draw[arrows=->] (12r)--(-12r);
\draw[arrows=->] (-12r)--(-32r);

\draw[arrows=->] (33r)--(13r);
\draw[arrows=->] (13r)--(-13r);
\draw[arrows=->] (-13r)--(-33r);

\draw[arrows=->] (34r)--(14r);
\draw[arrows=->] (14r)--(-14r);
\draw[arrows=->] (-14r)--(-34r);
\draw[arrows=->] (-34r)--(-54r);

\draw[arrows=->] (55r)--(35r);
\draw[arrows=->] (35r)--(15r);
\draw[arrows=->] (15r)--(-15r);
\draw[arrows=->] (-15r)--(-35r);
\draw[arrows=->] (-35r)--(-55r);

\node (h) at (-2,-2) {$h_j=$};
\node (h) at (0,-2) {$4$};
\node (h) at (2,-2) {$2$};
\node (h) at (4,-2) {$0$};
\node (h) at (6,-2) {$-2$};
\node (h) at (8,-2) {$-4$};
\node (h) at (10,-2) {$-6$};
\node (h) at (14,-2) {$5$};
\node (h) at (16,-2) {$3$};
\node (h) at (18,-2) {$1$};
\node (h) at (20,-2) {$-1$};
\node (h) at (22,-2) {$-3$};
\node (h) at (24,-2) {$-5$};
\node (h) at (26,-2) {$-7$};

\textcolor{red}{
\node (h) at (-2,-3) {$h_j=$};
\node (h) at (0,-3) {$5$};
\node (h) at (2,-3) {$3$};
\node (h) at (4,-3) {$1$};
\node (h) at (6,-3) {$-1$};
\node (h) at (8,-3) {$-3$};
\node (h) at (10,-3) {$-5$};
\node (h) at (14,-3) {$6$};
\node (h) at (16,-3) {$4$};
\node (h) at (18,-3) {$2$};
\node (h) at (20,-3) {$0$};
\node (h) at (22,-3) {$-2$};
\node (h) at (24,-3) {$-4$};
\node (h) at (26,-3) {$-6$};
}
\end{tikzpicture}
\caption{$\tilde{H}^0(G\times_BV)$ and ${H}^0(G\times_BV)$.}
\label{fig: H0 and tildeH0}
\end{figure}
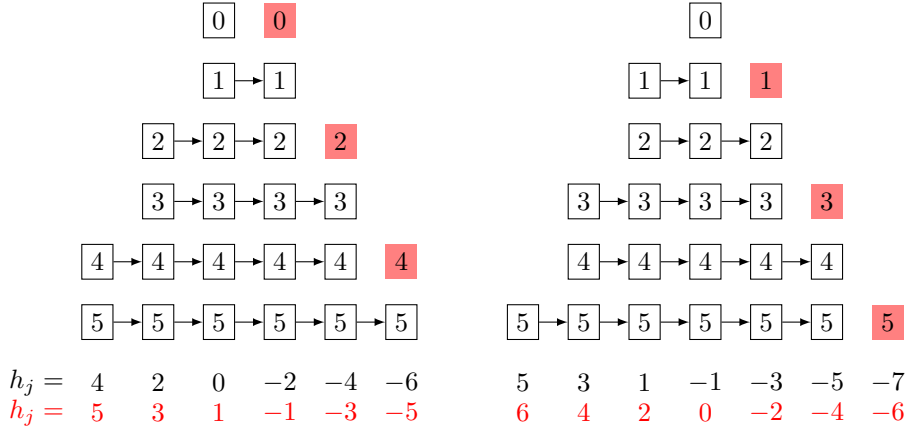

Let us assume that there exists $B$-modules
\begin{align*}
\underbrace{\tilde{H}^0(G\times_B\cdots\tilde{H}^0(G\times_B}_{m}\tilde{V}_{\hat{\lambda};\lambda_1,\dots,\lambda_N})\cdots)\quad (0\leq m\leq N-1)
\end{align*} 
with the conformal grading 
such that
\begin{enumerate}
\item 
For a dominant integral weight $\beta\in P_+$ and $\sigma\in W$, we have
\begin{align} 
\begin{split}\label{condition 1}
&\ch_q\underbrace{\tilde{H}^0(G\times_B\cdots\tilde{H}^0(G\times_B}_{m\geq 0} \tilde{V}_{\hat{\lambda};\cdots,\lambda_{m+1},\cdots})\cdots)^{h=\sigma\circ \beta}\\
=&
\ch_q\underbrace{\tilde{H}^0(G\times_B\cdots\tilde{H}^0(G\times_B}_{m\geq 0}\tilde{V}_{{\hat{\lambda}-\epsilon_{\lambda_{m+1}}(\sigma)};\cdots,\sigma\ast\lambda_{m+1},\cdots})\cdots)^{h=\beta-\epsilon_{\lambda_{m+1}}(\sigma)}.
\end{split}
\end{align}
\item 
For $\beta_m\in P_+$ such that $\beta_m-(N-m-1)\rho\in P_+$, we have
\begin{align}
&\ch_q(\tilde{V}_{{\hat\lambda};\lambda_1,\dots,\lambda_N})^{h=\beta_0}
=\ch_q({V}_{{\hat\lambda};\lambda_1,\dots,\lambda_N})^{h=\beta_0},\label{condition 2}\\
\begin{split}\label{condition 3}
&\ch_q\underbrace{\tilde{H}^0(G\times_B\tilde{H}^0(G\times_B\cdots\tilde{H}^0(G\times_B}_{m\geq 2} \tilde{V}_{{\hat\lambda};\cdots,\lambda_{m},\cdots})\cdots)^{h=\beta_m}\\
=&
\ch_q\underbrace{{H}^0(G\times_B\tilde{H}^0(G\times_B\cdots\tilde{H}^0(G\times_B}_{m\geq 2}\tilde{V}_{{\hat\lambda+\rho};\cdots,\lambda_m,\cdots})\cdots)^{h=\beta_m+\rho},\\
&\ch_q\tilde{H}^0(G\times_B\tilde{V}_{\hat\lambda;\lambda_1,\dots})^{h=\beta_1}
=\ch_q{H}^0(G\times_B\tilde{V}_{\hat\lambda+\rho;w_0\ast\lambda_1,\dots})^{h=\beta_1+\rho}.
\end{split}
\end{align}
\item 
We have the cohomology vanishings
\begin{align}\label{condition 4}
H^1(G\times_B\underbrace{\tilde{H}^0(G\times_B\cdots\tilde{H}^0(G\times_B}_{m\geq 0}\tilde{V}_{{\hat\lambda};\lambda_1,\dots,\lambda_N})\cdots)=0.
\end{align}
In particular, we can use the Atiyah--Bott fixed point formula \cite{AB} 
\begin{align*}
\sum_{a\geq 0}(-1)^n\ch_q H^a(G\times_BV)
=
\sum_{\beta\in P_+}
\dim L(\beta)
\sum_{\sigma\in W}
(-1)^{l(\sigma)}\ch_qV^{h=\sigma\circ\beta}
\end{align*}
to compute the characters of $0$-th sheaf cohomology $H^0(G\times_B-)$.
\end{enumerate}
Our goal is to compute the character of $H^0(G\times_B\underbrace{\tilde{H}^0(G\times_B\cdots\tilde{H}^0(G\times_B}_{N-1}\tilde{V}_{\vec{\lambda}})\cdots)$.
Conditions \eqref{condition 1} and \eqref{condition 4} are understood from the viewpoint of Felder exactness.
Namely, it is expected that for each $0\leq m\leq N-1$, there exists a Felder exact $\{V_{\lambda_m},Q^{[\lambda_m]}_j\}_{\lambda_m\in\Lambda_m,j\in I}$ such that
\begin{align*}
V_{\lambda_m}=\underbrace{\tilde{H}^0(G\times_B\cdots\tilde{H}^0(G\times_B}_{m}\tilde{V}_{{\hat\lambda};\dots,\lambda_m,\dots})\cdots).
\end{align*}
The meaning of \eqref{condition 2} is clear. 
Condition \eqref{condition 3} is illustrated as
Figure \ref{fig: H0 and tildeH0}.
Here,  
{\protect\tikz[baseline=(T.base)]\protect\node[fill=red!50](T){$k$};}
represents the difference between $\tilde{H}^0(G\times_BV)$ and ${H}^0(G\times_BV)$.
In ${H}^0(G\times_BV)$,
$(
\begin{tikzpicture}
[baseline=(00r.base)]
\node[draw] (00r) at (1,0) {$k+1$};
\node[draw] (01l) at (0,0) {$k$};
\end{tikzpicture}
)$
(or
$(
\begin{tikzpicture}
[baseline=(00r.base)]
\node[draw] (00r) at (1,0) {$k+1$};
\node[fill=red!50] (01l) at (0,0) {$k$};
\end{tikzpicture}
)$
)
is considered as a family of weight vectors of simple $\mathrm{SL}_2^j$-module with highest weight $k\varpi$.
Note that adding 
\begin{tikzpicture}
[baseline=(01l.base)]
\node[fill=red!50] (01l) at (0,0) {$k$};
\end{tikzpicture}
does not change the subspace $h_j\geq 0$, and that Cartan weight shifts by $\rho$ between $\tilde{H}^0(G\times_BV)$ and ${H}^0(G\times_BV)$.
For the reason of $\lambda_1\mapsto w_0\ast\lambda_1$, see \cite{KSug}.

\subsection*{Example 1}
($(1,p)$-log VOA, see \cite{FT,Sug}).
Let us consider the case $N=1$, and denote $\vec{\lambda}$ simply as $\lambda$.
In this case, $V_{\hat\lambda;\lambda}$ has the $B$-action by
\begin{align*}
f_i=\int e^{\sqrt{p}\alpha_i}(z)dz,
\quad
h=\lceil-\tfrac{1}{\sqrt{p}}{\alpha_i}_{(0)}\rceil.
\end{align*}
The $B$-action commutes with the short screening operators
$Q^{[\lambda]}_{j}\colon V_{\lambda}\rightarrow V_{\sigma_j\ast\lambda}$, and $\{V_{\lambda},Q^{[\lambda]}_{j}\}_{\lambda\in\Lambda,j\in I}$ is Felder exact.
The condition \eqref{condition 3} does not appear in this case and \eqref{condition 1}, \eqref{condition 2} are clear.
As explained above, \eqref{condition 4} was proved for special cases $(-p\lambda+\rho,\theta)\leq p$ in \cite{Sug}.

\subsection*{Example 2}
($(p_1,p_2)$-log VOA, see \cite[Section 6.2]{KSug}.)
Let us consider the case $N=2$.
The representation theory of $(p_1,p_2)$-log VOA is studied algebraically for the case of $\g=\sll_2$ \cite{FGST,TW}, but not for general $\g$ or geometrically.
Let us consider the case $\g=\sll_2$.
Unlike the case of $N=1$, $V_{\vec{\lambda}}$ cannot consists a Felder exact sequence, because of the shape of socle sequence \cite[Figure 3]{KSug}. 
However, the image of the short screening operator $Q^{[\lambda_1]}$ has the $B$-module structure definend by the Frobenius homomorphism \cite{TW}, and combining another short screening operator $Q^{[\lambda_2]}$ defines a Felder exact set.
Hence we can take $\tilde{H}^0(G\times_BV_{\vec{\lambda}})$ as $\im Q^{[\lambda_1]}$.
However, the structure of $\tilde{V}_{\vec{\lambda}}$ is not well understood.
In \cite{KSug}, we gave a makeshift definition of  $\tilde{V}_{\vec{\lambda}}$ that works for the computation of character below, but from a viewpoint of free field realization of log VOAs, $\tilde{V}_{\vec{\lambda}}$ should be a direct sum of Fock modules. 
In that case, however, the $B$-action would be complicated and not be a Virasoro homomorphism, although the conformal grading be preserved.

\begin{remark}
There is another reason to believe that our B-actions might not be Virasoro homomorphisms in general.
For simplicity, let us consider the case $\g=\sll_2$.
As we will see below, our character computation is to divide $\tilde{V}_{\vec{\lambda}}$ (or  Fock modules) into $2^{N}$ kinds of components and reduce them by half each time $H^0(G\times_B-)$ is taken (see Figure \ref{fig: flow of computation}).
In fact, when $N=1,2$, the socle sequence of $\tilde{V}_{\vec{\lambda}}$ consists of $2$ and $4$ types of irreducible Virasoro modules, respectively, and the method above results in a direct sum of one type of irreducible Virasoro module at the end.
However, the decomposition of Fock modules by irreducible Virasoro modules is the finest for $N=2$ (i.e., minimal model).
Therefore, a finer decomposition of $\tilde{V}_{\vec{\lambda}}$ than Virasoro is expected in the case of $N\geq 3$.
\end{remark}

\subsection{Character formula}\label{character formula}
\begin{figure}[htbp]
\centering
\begin{tikzpicture}[scale=.4]

\node (z) at (-2,16) {$\tilde{V}_{\Vec{\lambda}}=$};

\node[draw]
(-43r) at (8,14) {$3$};
\node[draw]
(-44r) at (8,12) {$4$};

\node[draw]
(-21r) at (6,18) {$1$};
\node[draw]
(-22r) at (6,16) {$2$};
\node[draw]
(-23r) at (6,14) {$3$};
\node[draw]
(-24r) at (6,12) {$4$};

\node[draw]
(00r) at (4,20) {$0$};
\node[draw]
(01r) at (4,18) {$1$};
\node[draw]
(02r) at (4,16) {$2$};
\node[draw]
(03r) at (4,14) {$3$};
\node[draw]
(04r) at (4,12) {$4$};

\node[draw]
(22r) at (2,16) {$2$};
\node[draw]
(23r) at (2,14) {$3$};
\node[draw]
(24r) at (2,12) {$4$};

\node[draw]
(44r) at (0,12) {$4$};

\draw[arrows=->] (01r)--(-21r);

\draw[arrows=->] (22r)--(02r);
\draw[arrows=->] (02r)--(-22r);

\draw[arrows=->] (23r)--(03r);
\draw[arrows=->] (03r)--(-23r);
\draw[arrows=->] (-23r)--(-43r);

\draw[arrows=->] (44r)--(24r);
\draw[arrows=->] (24r)--(04r);
\draw[arrows=->] (04r)--(-24r);
\draw[arrows=->] (-24r)--(-44r);

\node[draw]
(-44l) at (20,12) {$4$};
\node[draw]
(-22l) at (18,16) {$2$};
\node[draw]
(-24l) at (18,12) {$4$};
\node[draw]
(00l) at (16,20) {$0$};
\node[draw]
(02l) at (16,16) {$2$};
\node[draw]
(04l) at (16,12) {$4$};
\node[draw]
(22l) at (14,16) {$2$};
\node[draw]
(24l) at (14,12) {$4$};
\node[draw]
(44l) at (12,12) {$4$};
\draw[arrows=->] (22l)--(02l);
\draw[arrows=->] (02l)--(-22l);
\draw[arrows=->] (44l)--(24l);
\draw[arrows=->] (24l)--(04l);
\draw[arrows=->] (04l)--(-24l);
\draw[arrows=->] (-24l)--(-44l);

\draw[arrows=->] (9,16)--(11,16);
\node (a) at (10,17) {$H^0(G\times_B-)$}; 

\node (h) at (-2,10) {$h_j=$};
\node (h) at (0,10) {$4$};
\node (h) at (2,10) {$2$};
\node (h) at (4,10) {$0$};
\node (h) at (6,10) {$-2$};
\node (h) at (8,10) {$-4$};
\node (h) at (12,10) {$4$};
\node (h) at (14,10) {$2$};
\node (h) at (16,10) {$0$};
\node (h) at (18,10) {$-2$};
\node (h) at (20,10) {$-4$};

\node (b') at (-2,5) {$\text{Figure \ref{fig: H0 and tildeH0}}$};
\node (b) at (-2,4) {$=$};

\node[draw] (44) at (0,0) {$4$};
\node[draw] (24) at (2,0) {$4$};
\node[draw] (04) at (4,0) {$4$};
\node[draw] (-24) at (6,0) {$4$};
\node[draw] (-44) at (8,0) {$4$};

\node[draw] (43) at (0,1.2) {$3$};
\node[draw] (23) at (2,1.2) {$3$};
\node[draw] (03) at (4,1.2) {$3$};
\node[draw] (-23) at (6,1.2) {$3$};
\node[fill=red!50] (-43) at (8,1.2) {$3$};

\node (443) at (0,0.6) {$\ $};
\node (243) at (2,0.6) {$\ $};
\node (043) at (4,0.6) {$\ $};
\node (-243) at (6,0.6) {$\ $};
\node (-443) at (8,0.6) {$\ $};

\draw[arrows=->] (443)--(243);
\draw[arrows=->] (243)--(043);
\draw[arrows=->] (043)--(-243);
\draw[arrows=->] (-243)--(-443);

\node[draw] (22) at (2,3) {$2$};
\node[draw] (02) at (4,3) {$2$};
\node[draw] (-22) at (6,3) {$2$};

\node[draw] (21) at (2,4.2) {$1$};
\node[draw] (01) at (4,4.2) {$1$};
\node[fill=red!50] (-21) at (6,4.2) {$1$};

\node (221) at (2,3.6) {$\ $};
\node (021) at (4,3.6) {$\ $};
\node (-221) at (6,3.6) {$\ $};

\draw[arrows=->] (221)--(021);
\draw[arrows=->] (021)--(-221);

\node[draw] (00) at (4,6) {$0$};

\draw[stealth-stealth] (7,4)--(9,4);
\node (c) at (8,5) {$\eqref{condition 3}$};

\node (y) at (12.5,4) {$\tilde{H}^0(G\times_BV_{\vec{\lambda'}})=$};

\node[draw] (44') at (15,0) {$4$};
\node[draw] (24') at (17,0) {$4$};
\node[draw] (04') at (19,0) {$4$};
\node[draw] (-24') at (21,0) {$4$};
\node[draw] (-44') at (23,0) {$4$};

\node[draw] (43') at (15,2) {$3$};
\node[draw] (23') at (17,2) {$3$};
\node[draw] (03') at (19,2) {$3$};
\node[draw] (-23') at (21,2) {$3$};

\draw[arrows=->] (44')--(24');
\draw[arrows=->] (24')--(04');
\draw[arrows=->] (04')--(-24');
\draw[arrows=->] (-24')--(-44');

\draw[arrows=->] (43')--(23');
\draw[arrows=->] (23')--(03');
\draw[arrows=->] (03')--(-23');

\node[draw] (22') at (17,4) {$2$};
\node[draw] (02') at (19,4) {$2$};
\node[draw] (-22') at (21,4) {$2$};

\node[draw] (21') at (17,6) {$1$};
\node[draw] (01') at (19,6) {$1$};

\draw[arrows=->] (22')--(02');
\draw[arrows=->] (02')--(-22');
\draw[arrows=->] (21')--(01');

\node[draw] (00') at (19,8) {$0$};

\node (h) at (-2,-2) {$h_j=$};
\node (h) at (0,-2) {$4$};
\node (h) at (2,-2) {$2$};
\node (h) at (4,-2) {$0$};
\node (h) at (6,-2) {$-2$};
\node (h) at (8,-2) {$-4$};
\node (h) at (15,-2) {$3$};
\node (h) at (17,-2) {$1$};
\node (h) at (19,-2) {$-1$};
\node (h) at (21,-2) {$-3$};
\node (h) at (23,-2) {$-5$};

\draw[arrows=->] (24,4)--(26,4);
\node (d) at (25,5) {$H^0(G\times_B-)$}; 

\node (e) at (27,4) {$\cdots$};
\end{tikzpicture}
\caption{Flow of computation for $\ch_qH^0(G\times_B\tilde{H}^0(G\times\cdots V_{\vec{\lambda}})\cdots)$.}
\label{fig: flow of computation}
\end{figure}
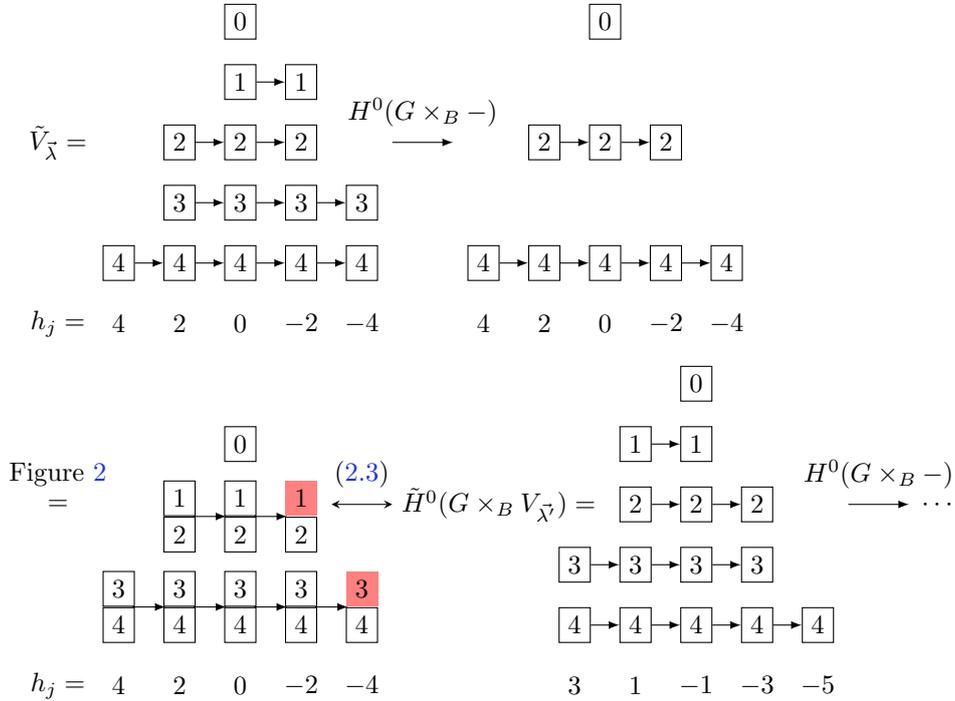

Using \eqref{condition 1}$\sim$\eqref{condition 4}, we can compute the character as
\begin{align*}\label{the character of N-fibered case for general case}
&\eta(q)^{\mathrm{rank}\g}\ch_qH^0(G\times_B\underbrace{\tilde{H}^0(G\times_B\cdots\tilde{H}^0(G\times_B}_{N-1} V_{{\hat\lambda};\lambda_1,\dots,\lambda_N})\cdots)^{h=\hat\lambda}\\
=
&\sum_{\beta_N\in P_+}\sum_{w_N\in W}(-1)^{l(w_N)}m_{\beta_N,\hat\lambda}\sum_{\beta_{N-1}\in P_+}m_{\beta_{N-1},\beta_N-\epsilon_{\lambda_N}(w_N)+\rho}\sum_{w_{N-1}\in W}(-1)^{l(w_{N-1})}\\
&\sum_{\beta_{N-2}\in P_+}m_{\beta_{N-2},\beta_{N-1}-\epsilon_{\lambda_{N-1}}(w_{N-1})+\rho}\sum_{w_{N-2}\in W}(-1)^{l(w_{N-2})}\\
&\cdots\sum_{\beta_{1}\in P_+}m_{\beta_1,\beta_2-\epsilon_{\lambda_{2}}(w_{2})+\rho}\sum_{w_{1}\in W}(-1)^{l(w_{1})}\ch_qV_{w_1\ast w_0\ast \lambda_1,w_2\ast\lambda_2,\dots,w_N\ast \lambda_N}^{h=\beta_1-\epsilon_{w_0\ast\lambda_{1}}(w_1)}\\
=
&\sum_{\gamma_1,\dots,\gamma_{N}\geq 0}p(\gamma_1)\cdots p(\gamma_{N})
\sum_{w_1,\dots,w_N\in W}(-1)^{l(w_1)+\cdots+l(w_N)}\\
&\ch_qV_{w_1\ast w_0\ast \lambda_1,w_2\ast\lambda_2,\dots,w_N\ast \lambda_N}^{h=\gamma_1+\cdots+\gamma_N+\hat\lambda+(N-1)\rho-\epsilon_{w_0\ast\lambda_1}(w_1)-\sum_{i=2}^{N}\epsilon_{\lambda_i}(w_i)},
\end{align*}
where $p(-)$ and $m_{-,-}$ are Kostant's partition function and multiplicity, respectively.
In the case of $\g=\sll_2$, this computation is illustrated as Figure \ref{fig: flow of computation}.
An appropriate linear combination of them would give \cite[(4.32)]{CFGHP}.

\subsection*{Example 3}
Let us restrict ourselves to the case of $\g=\sll_2$.
Then we have
\begin{align*}
&\eta(q)\ch_q H^0(G\times_B\tilde{H}^0(G\times_B\cdots V_{\vec{\lambda}})\cdots)^{h=\hat\lambda}\\
=&\sum_{n\geq 0}\binom{n+N-1}{N-1}\sum_{\epsilon_1,\dots,\epsilon_N\in\{\pm 1\}}\epsilon_1\cdots\epsilon_N
q^{\tfrac{p}{4}(2n+N+s-1-\sum_{i=1}^N\tfrac{\epsilon_i r_i}{p_i})^2},
\end{align*}
where $\hat\lambda=(s-1)\varpi$ for $s\in\{1,2\}$.
An appropriate linear combination of them would give the homological block of the corresponding $(N+2)$-fibered Seifert manifold \cite[Definition 2.1]{MTer}.

\subsection*{Example 4}
Let us restrict ourselves to the case of $N=2$.
As noted above, the $(p_1,p_2)$-log VOA and its modules associated with general $\g$ have not yet been constructed, but their characters can be inferred as
\begin{align*}
\sum_{\gamma_1,\gamma_2\geq 0}p(\gamma_1)p(\gamma_2)
\sum_{w_1,w_2\in W}(-1)^{l(w_1)+l(w_2)}
\ch_qV_{w_1\ast w_0\ast\lambda_1,w_2\ast\lambda_2}^{h=\gamma_1+\gamma_2+\hat\lambda+\rho-\epsilon_{w_0\ast\lambda_1}(w_1)-\epsilon_{\lambda_2}(w_2)}.
\end{align*}
We hope to generalize \cite[\S5.3]{KSug} to a relationship between this formula with $\g$-quantum invariant for torus link $T_{p_1\rk,p_2\rk}$ as a generalization of \cite{K}.

\end{document}